\documentclass[11pt]{amsart}
\usepackage[utf8]{inputenc}
\usepackage{amsfonts,latexsym,amssymb,amsmath,amscd,amsthm}
\usepackage[dvips]{graphicx}
\usepackage{enumerate}
\usepackage{amscd}

\newtheorem{teorema}{Theorem}
\newtheorem{propo}[teorema]{Proposition}

\newtheorem{coro}[teorema]{Corollary}

\theoremstyle{definition}

\theoremstyle{remark}
\newtheorem{nota}[teorema]{Remark}

\makeatletter          
\def\@@and{y}
\renewcommand{\andify}{%
\nxandlist{\unskip, }{\unskip{} \@@and~}{\unskip \@@and~}}
\def\and{\unskip{ }\@@and{ }\ignorespaces}
\makeatother           

\newcommand{\be}{\begin{enumerate}}
\newcommand{\ee}{\end{enumerate}}

\newcommand{\bi}{\begin{itemize}}
\newcommand{\ei}{\end{itemize}}

\newcommand{\norm}[1]{\left | #1 \right |}

\newcommand{\RR}{{\mathbb R}}
\newcommand{\ZZ}{{\mathbb Z}}
\newcommand{\NN}{{\mathbb N}}

\newcommand{\CC}{{\mathbb C}}

\newcommand{\bo}{{\mathbf 0}}

\begin{document}

\title{Entire functions polynomially bounded in several variables}

\author{Jorge Mozo-Fern\'{a}ndez}
\address[Jorge Mozo Fern\'{a}ndez]{Dpto. Álgebra, Análisis Matemático, Geometría y Topología \\
Facultad de Ciencias \\ Campus Miguel Delibes \\
Universidad de Valladolid\\
Paseo de Bel\'{e}n, 7 \\
47011 Valladolid\\
Spain}
\email{jorge.mozo@uva.es}
\thanks{Corresponding author: Jorge Mozo-Fernández (e-mail: jorge.mozo@uva.es)}.
\thanks{The author is partially supported by the Ministerio de Ciencia e Innovación from Spain, under the Project ``\'{A}lgebra y Geometr\'{\i}a en Din'{a}mica Real y Compleja III" (Ref.: MTM2013-46337-C2-1-P)}




\date{\today}


\begin{abstract}
In this paper we show that if an entire function $f(z_1,z_2)$ of two (or more) complex variables verifies $\norm{f(z_1,z_2)}\leq K(\norm{P(z_1,z_2)})$, where $P(z_1,z_2)$ is a polynomial that is not a power in $\CC[[z_1,z_2]]$, and $K$ is any positive-valued real function, then $f(z_1,z_2)$ can be written as a holomorphic function of $P$.
\end{abstract}

\maketitle

\section{Introduction} \label{introduccion}

Entire functions constitute a very interesting subject in the
theory of holomorphic  functions. They have lots of properties
concerning their growth and zeros, that make them behave,
sometimes, in a rather curious way, specially if you put some restrictions
on them. For instance, one of the more striking results that a student of a
basic course on complex analysis learns is Liouville's Theorem,
that asserts that a bounded entire function must be constant.

Another interesting result is that an entire function $f(z)$ that
verifies bounds $\norm{f(z)}\leq C\cdot \norm{z}^N$ is
necessarily a polynomial of degree at most $N$, result that
follows easily from Cauchy Theorem and Cauchy estimates on the coefficients of the series expansion of the function $f(z)$.

This result may be generalized to several variables. For
instance, exercise E.4c in L. Kaup and B. Kaup textbook \cite{KK}
demands to prove that  an entire function $f(z_1,z_2,\ldots ,
z_n)$ that verifies $\norm{f(z_1,z_2,\ldots , z_n)}\leq
\gamma\cdot \norm{z_1^{\nu_1}z_2^{\nu_2}\cdots z_n^{\nu_n}} $ is
a polynomial of degree at most $\norm{\nu}:= \nu_1+\nu_2+\cdots +
\nu_n$. The proof is similar to the one variable case, and
moreover you can see that, in the power series expansion of $f$,  the terms 
$z_1^{r_1}z_2^{r_2}\cdots z_n^{r_n}$ appearing explicitely  must verify  $r_i\leq
\nu_i$, for each $i$, $1\leq i\leq n$.

A closer look to this result shows another properties. For
instance, $f(\bo )=0$, i.e., $f$ has no constant term. If
$\nu_i>0$, $f(z_1,\ldots ,z_{i-1},0,z_{i+1},\ldots , z_n)\equiv
0$, so every term of $f$ contains $z_i$, and you can divide
$$
\norm{\frac{1}{z_i} f(z_1,z_2,\ldots , z_n)} \leq \gamma \cdot
| z_1^{\nu_1}\cdots z_i^{\nu_1-1}\cdots z_n^{\nu_n} |.
$$
Iterating this, one can conclude that every term of $f$ must
contain the monomial $z_1^{\nu_1}z_2^{\nu_2}\cdots z_n^{\nu_n}$,
and so,
$f(z_1,z_2,\ldots , z_n)= a z_1^{\nu_1}z_2^{\nu_2}\cdots
z_n^{\nu_n}$. In order to arrive to this conclusion it is necessary to assume that the bound $\norm{f(z_1,z_2,\ldots , z_n)}\leq
\gamma\cdot \norm{z_1^{\nu_1}z_2^{\nu_2}\cdots z_n^{\nu_n}} $ is verified \textit{everywhere} in $\CC^n$, not only for big enough values of $|(z_1,z_2,\ldots , z_n)|$.

The above bound on $f$ can be  replaced by another bounds and one
can still obtain similar results. Let us illustrate this  in  a simple
case: suppose you have a bound
$$
\norm{f(z_1,z_2,\ldots , z_n)}\leq K(\norm{z_1^n}).
$$
Here, $f$ is an entire function, and $K:\RR_{\geq 0}\rightarrow
\RR_{\geq 0}$ is any function. If you fix $z_1=z_{10}$, you
obtain that  $f(z_{10},z_2,\ldots , z_n)$ is bounded, hence
constant. So, $f(z_1,z_2,\ldots , z_n)=h(z_1)$, $h$ entire. Note
that here, the exponent $n$ is not relevant, as one can choose
$K'=K\circ K_n$, where $K_n(x)=x^n$, and so,
$K'(\norm{z_1})=K(\norm{z_1^n})$. In general, this implies that it is not a restriction to impose that the polynomial that appears on the right hand side of the inequality is not a power of another power series.

The objective of this paper is to generalize previous results. We
will focus in two variable case, even if most of the results
remain true in an arbitrary number of variables. The main result
of the paper is the following:
\begin{teorema} \label{principal}
Let $f(z_1,z_2)$ be an entire function, $P(z_1,z_2)$ a polynomial
that is not a power considered as a power series, and $K:\RR_{\geq 0}\rightarrow \RR_{\geq 0}$
any function. If for every $(z_1,z_2)\in \CC^2$ there is a bound
$$
\norm{f(z_1,z_2)}\leq K(\norm{P(z_1,z_2)}),
$$
then $f$ is a holomorphic function of $P$, i.e., there exists an
entire function $h:\CC\rightarrow \CC$ such that  $f(z_1,z_2)=
h(P(z_1,z_2))$.
\end{teorema}

\begin{nota}
The hypothesis about not being a power implies in particular that $P(0,0)=0$.
\end{nota}

The levels of the polynomial $P(z_1,z_2) $ define an algebraic foliation ${\mathcal F}_P$ over the affine place $ \CC^2 $. So, in the language of the theory of holomorphic foliations, Theorem \ref{principal} says that if a function $f$ is bounded over every leaf of ${\mathcal F}_P$, then it is a holomorphic first integral of ${\mathcal F}_P$.

The proof we are going to present here uses techniques from algebraic geometry (Bertini's
Theorem and reduction of singularities) that we shall explain, in
order to guarantee that $f$ can be written as a function of $P$.
Holomorphy of $h$ can be shown using Inverse Function Theorem.

Replacing $P(z_1,z_2)$ by an arbitrary entire function (not a
polynomial), it seems that the result remains valid, but
different techniques should be used and we will not treat this
generalized case here.

There are another interesting results concerning entire functions and their behaviour. Let us mention, for instance, M. Suzuki's results \cite{Suzuki} concerning topological properties of bivariate complex polynomials, that can be extended to a certain class of entire functions, using results developed by T. Nishino \cite{Nishino1,Nishino2} about the surfaces defined by entire functions, relating them with their growth. Nevertheless, these results are of a different nature or ours.

The author would like to thank Sergio Carrillo (University of Valladolid) and Javier Rib\'{o}n (Universidade Federal Fluminense) for some interesting discussions on the subject.

\section{Proof of Theorem \ref{principal}}

Consider the hypothesis of Theorem \ref{principal}.
The first objective of this section will be to show that $f(z_1,z_2)$ is
constant over every leaf of the foliation over $\CC^2$ defined by
the levels of $P(z_1,z_2)$ and the hypothesis made on $P$
(i.e., $P(z_1,z_2)$ is not a proper power of any other polynomial)
imply that, in fact, $f(z_1,z_2)$ is constant over every level
set of $P(z_1,z_2)$.

For, consider the family of affine curves $P(z_1,z_2)=c$ in
$\CC^2$. These curves are not necessarily irreducible (and
hence, not necessarily connected). At this point, recall the first Bertini
Theorem, as stated in \cite{Shafarevich}.

\begin{teorema}[Bertini]
Consider a regular map $\sigma : X\rightarrow  Y$ between
irreducible algebraic varieties over a field $k$ of
characteristic 0. Let $\sigma^{\ast}:k (Y)\rightarrow k(X)$ be
the induced map over the function fields, and assume that $X$ is
irreducible over the algebraic closure $\overline{k(Y)}$ of
$k(Y)$. Then, $\sigma^{-1}(y)$  is generically irreducible, i.e.,
there exists  a non-empty Zariski open set $U\subseteq Y$ such
that, for $y\in U$, $\sigma^{-1}(y)$ is irreducible.
\end{teorema}

Let us apply this theorem to the following situation: $X$ will be
the algebraic variety over $\CC^3$ defined by the equation
$C-P(z_1,z_2)$ (with coordinates $(z_1,z_2,C)\in \CC^3$),
$Y=\CC$, and $\sigma: X\rightarrow Y$ is the projection over the
third coordinate. As $\CC (Y)=\CC (C)$, the algebraic closure
$\overline{\CC (C)}$ is contained in the field of Puiseux series,
i.e.,
$$
\overline{\CC (C)} \subseteq \bigcup_{n=1}^{\infty } \CC [[
C^{1/n}]][C^{-1}].
$$

If $C\in \CC$, $\sigma^{-1}(C)$ is the curve $P(z_1,z_2)=C$. In
order to apply Bertini's Theorem it is necessary to guarantee
that $X$ is irreducible over $\overline{\CC (C)}$. As $X$ is a
hypersurface defined by the equation $C-P(z_1,z_2)$, reducibility
of $X$ is equivalent to reducibility of the polynomial
$C-P(z_1,z_2)$ over $\CC (z_1,z_2)[ C^{1/n}]$ for some $n\in
\NN$, and this can happen only if the polynomial $C^n-P(z_1,z_2)$
is reducible.

\begin{propo}
Let $K$ be a field of characteristic 0, containing the roots of
unity, $a\in K$, and $P(X)=X^n-a\in K[X]$. If $P(X)$ is reducible
then there exist $b\in K$ and $m>1$ a divisor of $n$ such that
$b^m=a$. In other words, $a$ is a power in $K$.
\end{propo}

This result can be deduced from the theory of cyclic extensions, as stated in several texts of algebra and Galois theory. See, for instance, \cite{W} or the classical book of S. Lang \cite{Lang}. In order to be self-contained, let us provide a proof of this result.

\begin{proof}
Suppose that $P(X)=Q_1(X)Q_2(X)$ is a factorization of $P(X)$ in $K[X]$, where $Q_1$, $Q_2$ are monic
polynomials, and $Q_1$ is irreducible. Let $c$ be a $n$th root of  $a$ in
$\overline{K}$, algebraic closure of $K$, such that $c$ is also a root
of $Q_1(X)$. So, $Q_1$ turns out to be the irreducible polynomial
of $c$ over $K$.

The roots of $P(X)$ over $\overline{K}$ are $\{ c\cdot
\zeta_n^i;\ i=0,\ldots , n-1\}$, where $\zeta_n$ is a primitive
$n$th root of unity. Among them, $\{ c\cdot \zeta_n^{i_k};\
k=0,\ldots , r-1\}$ are the roots of $Q_1(X)$. The Galois
group $G$ of $Q_1(X)$ must be a group of order $r$, isomorphic to
a subgroup of $\ZZ /n\ZZ$. Indeed, if $\sigma \in G$, $\sigma$ is
determined by its value $\sigma (c)=c\cdot \zeta_n^{i(\sigma)}$,
and the identification $\sigma \mapsto i(\sigma)$ gives the
inclusion $G\hookrightarrow \ZZ /n\ZZ$. So, $\{ \zeta_n^{i_k};\
k=0,\ldots , r-1\}$ form a cyclic group, generated by $\zeta_n^m$,
with $mr=n$. Then, $c^r\in K$ and so, $a=c^n= (c^r)^m$.
\end{proof}

\begin{coro}
If $P(z_1,z_2)\in \CC [z_1,z_2]$ or $P(z_1,z_2)\in \CC [[z_1,z_2]]$ is not a power, then, for all
$n\in \NN$, $C^n-P(z_1,z_2)\in \CC [z_1,z_2,C]$ is irreducible.
In this situation, for every value of $c$ with at most finitely
many exceptions, the affine curve $P(z_1,z_2)=c$ is irreducible,
hence connected.
\end{coro}

Consider, now, the curve $P(z_1,z_2)=c$, for generic values of
$c$, and its projectivization ${\mathcal C}_c$ of equation $z_0^d \cdot
P\left(
\frac{z_1}{z_0}, \frac{z_2}{z_0}\right)-cz_0^d$, where $d$ is the
degree of $P$. By the hypothesis,
$$
\norm{f(z_1,z_2)}\leq K(\norm{P(z_1,z_2)})= K(\norm{c}),
$$
so, $f(z_1,z_2)$ is bounded over this curve. After reduction of
singularities, a compact connected Riemann surface $E_c$ is
obtained, and a map
$$
E_c\stackrel{\pi_c}{\longrightarrow } {\mathcal C}_c,
$$
that is a holomorphic diffeomorphism outside the singular locus
of ${\mathcal C}_c$, and such that each singular point has been
replaced by a finite number of points. So, $f\circ \pi_c$ is a
holomorphic bounded map defined on $E_c$ minus a finite number of
points. By the Riemann extension Theorem, $f\circ \pi_c$ admits a
unique holomorphic extension to $E_c$, that will be constant with
value $h(c)$. This function $h(c)$ can be extended by continuity
to the (finite in number) values of $c$ such that $P(z_1,z_2)=c$ is not
irreducible. So, it satisfies everywhere
$$
f(z_1,z_2)=h(P(z_1,z_2)),
$$
as desired.

It remains to show the holomorphy of the function $ h $  constructed above.
For, take $c$ such that
$P(z_1,z_2)=c$ is irreducible and $(z_{10}, z_{20})$ a regular
point of $P(z_1,z_2)=c$. Assume that $\dfrac{\partial
P}{\partial z_1}(z_{10}, z_{20})\neq 0$. The map $(z_1,z_2)
\rightarrow (P(z_1,z_2), z_2)$ is, then, a diffeomorphism in a
neighbourhood of $ (z_{10}, z_{20})$, and so, by Inverse Function
Theorem, the exists (locally) an inverse $F(c,z_2)= (F_1(c,z_2),
z_2 )$. Then
$$
f(F_1(c,z_2),z_2)= h(P(F_1(c,z_2), z_2))=h(c),
$$
and so, $h(c)$ is holomorphic.

\section{Conclusion and open questions}
In this note we have worked in a two variable setting, starting
with entire functions. The transition to more variables is not
difficult, and it is only a technical matter. If we replace
$P(z_1,z_2)$ by an entire function, not necessarily a polynomial,
the result seems to remain true, but this proof, based in some
results coming from algebraic geometry, is no longer valid.

Another question would be to relax the hypothesis about the
domain of definition of $f(z_1,z_2)$. In this context, when the
domain is a monomial sector, as defined in \cite{CDMS}, and
$P(z_1,z_2)$ is a monomial, similar results can be stated,
using different techniques related with the theory of monomial
asymptotic expansions,  developed in a 
paper of S. Carrillo and the author \cite{CM}.

\end{document}